\input amstex
\magnification=1200
\documentstyle{amsppt}
\NoRunningHeads
\NoBlackBoxes
\define\FA{\frak A}
\define\Fv{\frak v}
\define\Fg{\frak g}
\define\Fw{\frak w}
\define\Fp{\frak p}
\define\Fh{\frak h}
\define\Fb{\frak b}
\define\Gal{\operatorname{Gal}}
\define\End{\operatorname{End}}
\define\sothree{\operatorname{\frak s\frak o}(3)}
\define\sofour{\operatorname{\frak s\frak o}(4)}
\define\crc{\Bbb S^1}
\define\Vect{\operatorname{Vect}(\crc)}
\define\CVect{\Bbb C\Vect}
\define\sltwo{\operatorname{\frak s\frak l}(2,\Bbb C)}
\define\jltwo{\operatorname{\frak j\frak l}(2,\Bbb C)}
\define\slstwo{\operatorname{\frak s\frak l}^*(2,\Bbb C)}
\define\HS{\operatorname{\Cal H\Cal S}}
\define\SU{\operatorname{\bold S\bold U}}
\define\TC{\operatorname{\Cal T\Cal C}}
\define\RWsltwo{\operatorname{\Cal R\Cal W}(\sltwo)}
\define\RWhsltwo{\operatorname{\Cal R\Cal W}_3(\sltwo)}
\define\bnd{\Cal B}
\topmatter
\title Topics in hidden symmetries. V.
\endtitle
\author Denis V. Juriev
\endauthor
\affil\rm
``Thalassa Aitheria'' Research Center for Mathematical Physics
and Informatics,\linebreak
ul.Miklukho-Maklaya 20-180, Moscow 117437 Russia.\linebreak
E-mail: denis\@juriev.msk.ru\linebreak
\ \linebreak
\ \linebreak
November, 25, 1996.\linebreak
E-version: funct-an/9611003\linebreak
\ \linebreak
\ \linebreak
\endaffil
\abstract This note being devoted to some aspects of the inverse
problem of representation theory contains a new insight into it illustrated
by two topics. The attention is concentrated on the manner of representation
of abstract objects by the concrete ones as well as on the abstract objects
themselves.
\endabstract
\endtopmatter
\document

This paper being a continuation of the previous four parts [1] explicates
novel features of the general ideology presented in the review [2] (see also
the research article [3], where a wide interpretation of the 
inverse problem of representation theory was originally proposed). The 
attention is concentrated on the manner of concrete representations of 
abstract objects as well as on the least themselves.

The concrete lines of topics have their origins in the author's papers
[4,5:App.A], where some objects, which gave start to the constructions below,
appeared.

\head 1. Topic Nine: Lie composites and their representation. Composed
representations of Lie algebras
\endhead

This topic may be considered as a variation of a theme briefly discussed at the
end of the article [6] (the so--called isotopic composites and their
representations). Though the particular case considered below looses some
interesting combinatorial features of the general one, e.g. the
graph--representations, it is nevertheless rather interesting.

\definition{Definition 1}

{\bf A.} A linear space $\Fv$ is called {\it a Lie composite\/}
iff there are fixed its subspaces $\Fv_1,\ldots \Fv_n$ ($\dim\Fv_i>1$)
supplied by the compatible structures of Lie algebras. Com\-pa\-ti\-bi\-li\-ty 
means that
the structures of the Lie algebras induced in $\Fv_i\cap\Fv_j$ from $\Fv_i$ and
$\Fv_j$ are the same. The Lie composite is called {\it dense\/} iff
$\Fv_1\uplus\ldots\uplus\Fv_n=\Fv$ (here $\uplus$ denotes the sum of linear 
spaces). The Lie composite is called {\it connected\/}
iff for all $i$ and $j$ there exists a sequence $k_1,\ldots k_m$ ($k_1=i$,
$k_m=j$) such that $\Fv_{k_l}\cap\Fv_{k_{l+1}}\ne\varnothing$.

{\bf B.} {\it A representation\/} of the Lie composite $\Fv$ in the space $H$
is the linear mapping $T:\Fv\mapsto\End(H)$ such that $\left.T\right|_{\Fv_i}$
is a representation of the Lie algebra $\Fv_i$ for all $i$.

{\bf C.} Let $\Fg$ be a Lie algebra. A linear mapping $T:\Fg\mapsto\End(H)$
is called {\it the composed representation\/} of $\Fg$ in the linear space
$H$ iff there exists a set $\Fg_1,\ldots,\Fg_n$ of the Lie subalgebras of
$\Fg$, which form a dense connected composite and $T$ is its representation.
\enddefinition

Reducibility and irreducibility of representations of the Lie composites are
defined in the same manner as for Lie algebras. One may also formulate a
superanalog of the Definition 1. The set of representations of the fixed Lie 
composite is closed under the tensor product and, therefore, may be supplied by
the structure of {\it a tensor category}.

\example{Example 1 (The Octahedron Lie Composite)} Let us consider an
octahedron with the vertices
$A$, $B$, $C$, $D$, $E$, $F$, the edges $(AB)$, $(AC)$, $(AD)$, $(AE)$,
$(BC)$, $(BF)$, $(CD)$, $(CF)$, $(DE)$, $(DF)$, $(EF)$, and the faces
$(ABC)$, $(ACD)$, $(ADE)$, $(AEB)$, $(BCF)$, $(CDF)$, $(DEF)$, $(EBF)$.
Let $\Fv$ be a six--dimensional linear space with the basis labelled by the
vertices of the octahedron, $\Fv_1$, $\Fv_2$, $\Fv_3$, $\Fv_4$ be four
three--dimensional subspaces in $\Fv$ corresponded to the faces
$(ABC)$, $(ADE)$, $(CDF)$, $(EBF)$. All subspaces $\Fv_i$ are supplied by the
structures of the Lie algebras isomorphic to $\sothree$ (such structures are
compatible to the orientations on the faces). The pentuple
$(\Fv,\Fv_1,\Fv_2,\Fv_3,\Fv_4)$ is a dense connected Lie composite.

\proclaim{Proposition} Let $T$ be an arbitrary representation of the Lie
composite $(\Fv;\Fv_1,\Fv_2,\mathbreak\Fv_3,\Fv_4)$ in the finite--dimensional
linear space $H$, then $H$ admits a representation of the Lie algebra $\sofour$.
If $T$ is an irreducible representation then there exist the real numbers
$\lambda_A$, $\lambda_B$, $\lambda_C$, $\lambda_D$, $\lambda_E$, $\lambda_F$
such that the operators $T(A)-\lambda_A\boldkey 1$, $T(B)-\lambda_B\boldkey 1$,
$T(C)-\lambda_C\boldkey 1$, $T(D)-\lambda_D\boldkey 1$, $T(E)-\lambda_E\boldkey
1$, $T(F)-\lambda_F\boldkey 1$ form an irreducible representation of $\sofour$.
\endproclaim

\demo{Proof} First, note that the commutator of operators corresponded to
the opposite vertices commute with operators corresponded to other four
vertices. It commutes with all six operators because they may be expressed as
commutators of the least four operators. So the commutator of operators
corresponded to the opposite vertices belongs to the center of the Lie
algebra generated by the all six operators. Let us factorize this Lie algebra
$\Fg$ by the center. Such quotient is isomorphic to $\sofour$ (one uses the
fact that formulas for commutators of all six operators are known up to the
center of $\Fg$). The statement of the theorem is a consequence of this result
and the fact that any central extension of the semisimple Lie algebra is
trivial (i.e. may be splitted -- see f.e.[7])
\qed\enddemo
\endexample

The construction of the Octahedron Lie composite may be generalized on
the certain class of polyhedra. However, any analogs of the Proposition
are not known for such general case.

\example{Example 2 (The Witt composite)} Let $\Fw$ be the so--called Witt
algebra, which is a subalgebra of the complexification $\CVect$ of the Lie
algebra $\Vect$ of the smooth vector fields on a circle $\crc$ [8]. The Witt
algebra $\Fw$ consists of all polynomial vector fields and admits a basis
$e_k$ ($k\in\Bbb Z$) with commutation relations $[e_i,e_j]=(i-j)e_{i+j}$.

Let us consider two subalgebras $\Fp_{\pm}$ of $\Fw$ generated by $e_i$ with
$i\ge-1$ and $i\le1$; note that $\Fp_+\cap\Fp_-=\sltwo$. The triple
$(\Fw;\Fp_+,\Fp_-)$ is a dense connected Lie composite.

\proclaim{Theorem 1} The action of the Lie algebra $\sltwo$ in any Verma module
$V_h$ ($h$ is the highest weight) may be extended to the representation of the
Lie composite $(\Fw;\Fp_+,\Fp_-)$ and, hence, to the composed representation
of the Witt algebra $\Fw$.
\endproclaim

\demo{Proof} The proof of the theorem is based on the following Lemma.

\proclaim{Lemma [4:\S2.2;5:App.A3]} The generators of the Lie algebra $\sltwo$
in the Verma module $V_h$ may be included into the infinite family of tensor
operators (of the so--called $q_R$--conformal symmetries). If the Verma module
$V_h$ is realized in the space of polynomials of one variable $z$ and the
generators of the Lie algebra $\sltwo$ have the form
$$L_1=(\xi+2h)\partial_z,\quad L_0=\xi+h,\quad L_{-1}=z,$$
where $\xi=z\partial_z$ and $[L_i,L_j]=(i-j)l_{i+j}$, then the tensor operators
are of the form
$$L_k=(\xi+(k+1)h)\partial_z^k\quad (k\ge0),\quad
L_{-k}=z^k\frac{\xi+(k+1)h}{(\xi+2h)\ldots(\xi+2h+k-1)}\quad (k\ge 1).$$
\endproclaim

It follows from the explicit formulas for the tensor operators that
$[L_i,L_j]=(i-j)L_{i+j}$ for all $i,j\ge-1$ or $i,j\le1$
\qed\enddemo

\remark{Comment} The $\Fw$--representations of Theorem 1 in the Verma modules
over $\sltwo$ generate a tensor category of the category of all representations
of $\Fw$.
\endremark
\endexample

\remark{Remark 1} The construction of the Witt composite may be generalized
on the Riemann surfaces of higher genus in lines of I.M.Krichever and
S.P.Novikov [9].
\endremark

\remark{Remark 2} Generalizing the terminology of [1,2] one may say that
the tensor ope\-ra\-tors (of spin 2, i.e. the $q_R$--conformal symmetries
[4:\S2.2;5:App.A3]) in the Verma modules $V_h$ over the Lie algebra $\sltwo$
form the set of hidden symmetries, whose algebraic structure is one of the
Witt composite.
\endremark

Note that if the hidden symmetries realize a representation of the
Lie composite they should not be unpacked (a similar situation appears also
in the case of the isocommutator algebras of hidden symmetries and the
related Lie $\Fg$--bunches [1:To\-pic 3;2:\S2.1]).

\example{Example 3} Let $\Fw$ be the Witt algebra and $(\Fw;\Fp_{\pm})$ be the
Witt composite. Let us consider the abelian extension $\Fw^e$ of the Witt
algebra by the generators $f_i$ ($i\in\Bbb Z$) such that $[e_i,f_j]=jf_j$.
The subalgebras $\Fp_{\pm}$ of $\Fw$ may be extended to the subalgebras
$\Fp^e_{\pm}$ of $\Fw^e$ by the generators $f_i$, where $i\ge 0$ and $i\le 0$,
respectively. The triple $(\Fw^e;\Fp^e_{\pm})$ form the extended Witt composite.

\proclaim{Theorem 2} The representation of the Witt composite in any Verma 
module $V_h$ ($h$ is the highest weight) over the Lie algebra $\sltwo$ may be
extended to the representation of the Lie composite $(\Fw^e;\Fp^e_{\pm})$ and,
hence, the composed rep\-re\-sen\-tation of $\Fw$ in $V_h$ may be extended to 
the composed representation of $\Fw^e$.
\endproclaim

\demo{Proof} The additional generators $f_i$ are represented by the tensor
operators of spin 1, namely, $f_i\mapsto\partial_z^i$ ($i\ge 0$),
$f_{-i}\mapsto z^i\frac1{(\xi+2h)\ldots(\xi+2h+i-1)}$ ($i\ge 1$)
\qed\enddemo
\endexample

\remark{Comment} The $\Fw^e$--representations of Theorem 2 in the Verma modules 
over $\sltwo$ generate a tensor subcategory of the category of all 
representations of $\Fw^e$.
\endremark

\remark{Remark 3} It is very interesting to consider the composed
representations of the real semisimple Lie algebras $\Fg$, which unduce
representations of some natural sub\-al\-geb\-ras (for instance, of two opposite
maximal parabolic subalgebras or two opposite Borel subalgebras, perhaps plus
some $\sltwo$ imbed into $\Fg$, etc.).
\endremark

\head 2. Topic Ten: $\FA$--projective representations of Lie algebras
\endhead

The theme of this topics is $\FA$--projective representations ($\FA$ is an
associative algebra), which are certain generalizations of the ordinary
projective representations (see [9]).

\definition{Definition 2A} Let $\FA$ be an arbitrary associative algebra
represented in the linear space $H$ and $\Fg$ be a Lie algebra. The linear
mapping $T:\Fg\mapsto\End(H)$ is called {\it $\FA$--projective
representation\/} iff for all $X$ and $Y$ from $\Fg$ there exists an element
of $\FA$ represented by the operator $A_{XY}$ such that
$$[T(X),T(Y)]-T([X,Y])=A_{XY}.$$
\enddefinition

If $H$ is infinite dimensional the representation may be realized by the
unbounded operators.

\remark{Remark 4} The definition may be generalized on any anticommutative
algebras. In this situation it is closely related to the constructions of
representations of the anticommutative algebras $\jltwo$ and $\slstwo$ in
[10;11:\S2]. In general, it should be viewed in the context of the old ideas of
A.I.Maltsev on the representations of arbitrary nonassociative algebras [12].
The general anticommutative algebras and the constructions of their
$\FA$--projective representations is of an interest in the context of
the quasi--Hopf algebras, whose coalgebraic structure is not associative,
Lie (Jacobian) and co--Jacobian quasi--bialgebras and related structures [13]
(see also [14]).
\endremark

\remark{Example} If $(\Fg,\Fh)$ is the reductive pair then any representation
of $\Fg$ is an $\Cal U(\Fh)$--projective representation of the binary
anticommutative algebra $\Fp$ ($\Fg=\Fh\oplus\Fp$, the operation in $\Fp$ is
of the standard form: $[X,Y]_{\Fp}=\pi([X,Y])$, where $[\cdot,\cdot]$ is
the commutator in $\Fg$ and $\pi$ is the projector of $\Fg$ onto $\Fp$ along
$\Fh$ [15]).
\endremark

\remark{Remark 5} The standard projective representation appears if $\FA$
is one--dimensional algebra acting by scalar matrices.
\endremark

\remark{Remark 6} If $H$ is the Hilbert (or pre--Hilbert) space then one
may consider the algebra $\HS$ of all Hilbert--Schmidt operators as $\FA$.
In this case our construction is deeply related to ideas of A.I.Shtern
on the almost representations and to some aspects of the pseudodifferential
calculus [16] (see also [14]). One may also consider the algebra $\TC$ of all 
trace class operators and the algebra $\bnd$ of all bounded operators.
Any $\HS$--projective representation by unbounded operators is also a 
$\bnd$--projective representation.
\endremark

\definition{Definition 2B}
Let $\FA$ be an associative algebra with an involution $*$ symmetrically 
represented in the Hilbert space $H$. If $\Fg$ is a Lie algebra with an 
involution $*$ then its $\FA$--projective representation $T$ in the space $H$
is called {\it symmetric\/} iff for all elements $a$ from $\Fg$ $T(a^*)=T^*(a)$.
Let $\Fg$ be a $\Bbb Z$--graded Lie algebra ($\Fg=\oplus_{n\in\Bbb Z}\Fg_n$)
with an involution $*$ such that $\Fg^*_n=\Fg_{-n}$ and the involution is
identical on the subalgebra $\Fg_0$. Let us extend the $\Bbb Z$--grading and
the involution $*$ from $\Fg$ to the tensor algebra $\bold T^{\cdot}(\Fg)$.
The symmetric $\FA$--projective representation of $\Fg$ is called {\it
absolutely symmetric\/} iff for any element $a$ of $\bold T^{\cdot}(\Fg)$ 
such that $\deg(a)=0$ the equality $T(a)=T^*(a)$ holds (here the representation
$T$ of $\Fg$ in $H$ is extended to the mapping from $\bold T^{\cdot}(\Fg)$ to
$\End(H)$.
\enddefinition

Note that the Witt algebra $\Fw$ has the natural $\Bbb Z$--grading and
involution $*$.

\proclaim{Theorem 3} The composed representation of the Witt algebra $\Fw$
in the uni\-ta\-ri\-zable Verma module $V_h$ over the Lie algebra $\sltwo$ is
the absolutely symmetric $\HS$--projective representation of $\Fw$.
\endproclaim

The statement of the theorem follows from the explicit formulas for the
tensor operators of $q_R$--conformal symmetries in the Verma module over
the Lie algebra $\sltwo$.

It is very interesting to describe all highest weight composed representations
of the Witt composite, which are $\HS$--projective representations of the
Witt algebra or its central Gelfand--Fuchs extension (the Virasoro algebra) [8].

\remark{Remark 7} The Definition 2 may be generalized from the Lie algebras
to their numerous nonlinear analogs such as quantum groups, Sklyanin algebra,
Racah--Wigner algebras, mho--algebras, NWSO--algebras (see [1,2]), etc.
\endremark

It is very interesting to consider the composed representations of the
semisimple Lie algebras mentioned in Remark 3, which are the $\HS$--projective
representations of the algebras. The importance of such combination is
illustrated by the following example.

\remark{Example} Let $\Fb_{\pm}$ be two opposite Borel subalgebras of
$\sltwo$, then any rep\-re\-sen\-tation of $\Cal U_q(\sltwo)$ may be considered 
as a representation of the Lie composite $(\sltwo;\Fb_+,\Fb_-)$, however,
none of the infinite--dimensional $\Cal U_q(\sltwo)$--modules of the category
$\Cal O$ [17:\S9.2] defines the $\HS$--projective representation of $\sltwo$.
\endremark

\remark{Remark 8} The sets of the irreducible infinite--dimensional
$\HS$--projective highest weight representations for two nonlinear
$\operatorname{\frak s\frak l}_2$ [18] with commutation relations
$[e_0,e_{pm}]=\pm e_{\pm}$ and $[e_+,e_-]=R_{1,2}(e_0)$ coincides if and only
if for all $h\in\Bbb R$ $\sum_{j=0}^{\infty}|R_1(h+j)-R_2(h+j)|^2<\infty$.
\endremark

\definition{Definition 2C} The $\FA$--projective representation $T$ of the Lie
algebra $\Fg$ in the linear space $H$ will be called {\it almost absolutely 
closed\/} iff for any $n\ge 1$ and for any elements $X_0, X_1, X_2,\ldots 
X_{n+1}$ of $\Fg$ there exists an element $\varphi(X_0,X_1,X_2,\ldots X_{n+1})$ 
of $\Fg$ such that
$$[\ \ldots[[T(X_0),T(X_1)],T(X_2)],\,\ldots,T(X_{n+1})]\!\equiv\!
T(\varphi(X_0,X_1,X_2,\ldots X_{n+1}))\!\!\!\!\pmod{\FA},$$
here $\FA$ is considered as being mapped into $\End(H)$. The almost absolutely
closed $\FA$--projective representation $T$ of the Lie algebra $\Fg$ in the
linear space $H$ will be called {\it absolutely closed\/} iff
$\varphi(\cdot,\ldots,\cdot)\equiv 0$.
\enddefinition

The mappings $(X_0,X_1,X_2,\ldots X_{n+1})\mapsto\varphi(X_0,X_1,X_2,\ldots
X_{n+1})$ associated with any almost absolutely closed $\FA$--projective 
representation of the Lie algebra $\Fg$ define the higher brackets in the Lie 
algebra $\Fg$. The objects with higher brackets sys\-te\-ma\-ti\-cal\-ly 
appears in many branches of mathematics and mathematical physics (see, for
example, the book [14] and the article [19] among many others and numerous 
refs wherein).

\remark{Remark 9} The $\HS$--projective representations of the Witt algebra
in the Verma modules over $\sltwo$ are absolutely closed.
\endremark

Note that the extended Witt algebra $\Fw^e$ has the natural $\Bbb Z$--grading
and involution $*$.

\proclaim{Theorem 4} The composed representation of the extended Witt algebra 
$\Fw^e$ in the unitarizable Verma module $V_h$ over the Lie algebra $\sltwo$ 
is the absolutely closed and absolutely symmetric $\HS$--projective 
representation of $\Fw^e$.
\endproclaim

The statement follows from the explicit formulas for the tensor operators of
spins 1 and 2 in the Verma modules over the Lie algebra $\sltwo$.

\remark{Remark 10 (perhaps, very crucial)}
The $\HS$--projective (and, simultaneously, com\-posed) representations of the 
Lie algebras $\Fw$ and $\Fw^e$ from Theorems 3,4 realize generators of the
algebras as {\it infinite--dimensional hidden symmetries\/} of the Verma 
modules over the Lie algebra $\sltwo$. It seems that many other objects of
the representation theory of the reductive Lie algebras (Verma modules, their
models and skladens -- see e.g.[4]\footnote"*"{The originally Russian term 
``skladen'' was mesleadingly translated from Russian into English as 
``collection'' by a translator of the article [4], a possible correct 
translation is ``unfolding'' but the direct transliteration is 
preferable.\newline}, general constructive and Harish--Chandra modules, etc. 
-- see e.g.[20,17]) possess analogous infinite--dimensional hidden symmetries, 
which may constitute as familiar as novel abstract algebraic structures.
\endremark

\remark\nofrills{Problems:}
\roster
\item"--" What algebraic structure is represented by the tensor operators
$W_n$ ($n\in\Bbb Z$) of the spin 3 ($q_R$--$W_3$--symmetries) in the Verma
module $V_h$ over the Lie algebra $\sltwo$? The generators $W_i$
($i=-2,-1,0,1,2$) together with the $\sltwo$--generators form the Racah--Wigner
algebra $\RWsltwo$ for the Lie algebra $\sltwo$ [2:\S\S1.1,1.2;4:\S2.2]. The
conformally invariant analogue of the requested structure is the Zamolodchikov
$W_3$--algebra [21]. Some aspects of the problem were discussed in the
author's article [22:\S4].
\item"--" What algebraic structure is represented by the tensor operators
$U_n$ ($n\in\Bbb Z$) of the spin 4 in the Verma module $V_h$ over the Lie
algebra $\sltwo$? The generators $U_i$ ($i=-3,-2,-1,0,1,2,3$) together with
the $\sltwo$--generators form the mho--algebra $\mho(\sltwo,\pi_3)$ over
the Lie algebra $\sltwo$ [1:Topic 2;2:\S1.4]. The operators $U_i$
($i=-3,-2,-1,0,1,2,3$), $W_i$ ($i=-2,-1,0,1,2$) and the $\sltwo$--generators
form the higher Racah--Wigner algebra $\RWhsltwo$ [2:\S1.2] (see also
[23:\S1.2]).
\item"--" To describe the tensor category of representations of the Witt
composite generated by its representations in the Verma modules over
$\sltwo$ and by the Verma modules over the Gelfand--Fuchs central extension of 
the Witt algebra (Verma modules over the Virasoro algebra [24]).
\endroster
\endremark

\head Conclusions\endhead

Even the sketchy discussion of two topics above allows to state that the
actual richness and attractiveness of the inverse problem of representation
theory are based not only on a large scope of various interesting abstract
algebraic structures, which may be concretely represented and somehow
unravelled, but also on the diversity of the manners of representation,
which may produce very intriguing unexpected and nontrivial effects under
the intent look even in the simple and almost hackneyed situations.

\Refs
\roster
\item" [1]" Juriev D., Topics in hidden symmetries. I--IV. E-prints:
hep-th/9405050, q-alg/9610026, q-alg/9611003, q-alg/9611019.
\item" [2]" Juriev D.V., An excursus into the inverse problem of representation
theory [in Russian]. Report RCMPI-95/04 (August 1995) [e-version:
mp\_arc/96-477].
\item" [3]" Juriev D.V., Characteristics of pairs of operators, Lie hybrids,
Poisson brackets and non\-li\-near geometric algebra [in Russian].
Fundam.Prikl.Matem., to appear [e-version: funct-an/9411007].
\item" [4]" Juriev D.V. Complex projective geometry and quantum projective field
theory [in Russian]. Teor.Matem.Fiz. 101(3) (1994) 331-348 [English transl.:
Theor.Math.Phys. 101 (1994) 1387-1403].
\item" [5]" Juriev D.V., Belavkin--Kolokoltsov watch--dog effects in 
interactively controlled stochactic dynamical videosystems [in Russian]. 
Teor.Matem.Fiz. 106(2) (1996) 333-352 [English transl.: Theor.Math.Phys. 106 
(1996) 276-290].
\item" [6]" Juriev D.V., Topics in isotopic pairs and their representations. II.
A general supercase [in Russian]. Teor.Matem.Fiz., to appear.
\item" [7]" Guichardet A., Cohomologie des groupes topologiques et des
alg\`ebres de Lie. Paris, Cedic/Fernand Nathan, 1980.
\item" [8]" Fuchs D.B., Cohomology of the infinite dimensional Lie algebras
[in Russian]. Moscow, Nauka, 1984.
\item" [9]" Krichever I.M., Novikov S.P., Virasoro--type algebras, Riemann
surfaces and structures of the soliton theory [in Russian]. Funkts.analiz i
ego prilozh. 21(2) (1987) 46-63; Virasoro--type algebras, Riemann surfaces
and string in the Minkowsky space [in Russian]. ibid. 21(4) (1987) 47-61;
Virasoro--type algebras, the momentum--energy tensor and operator expansions
on the Riemann surfaces [in Russian]. ibid. 23(1) (1989) 1-14;
Virasoro--Gelfand--Fuchs algebras, Riemann surfaces, operator's theory of
closed strings. J.Geom. Phys. 5 (1988) 631-661 [reprinted in ``Geometry and
physics. Essays in honour of I.M.Gel\-fand'', Eds. S.Gindikin and I.M.Singer,
Pitagora Editrice, Bologna and Elsevier Sci.Publ., Amsterdam, 1991].
\item" [9]" Kirillov A.A., Elements of the representation theory. Springer,
1976.
\item"[10]" Juriev D., Noncommutative geometry, chiral anomaly in the quantum
projective ($\sltwo$--invariant) field theory and $\jltwo$--invariance.
J.Math.Phys. 33 (1992) 2819-2822, (E) 34 (1993) 1615.
\item"[11]" Juriev D.V., Quantum projective field theory: quantum--field
analogs of Euler--Arnold equations in projective $G$--hypermultiplets [in
Russian]. Teor.Matem.Fiz. 98(2) (1994) 220-240 [English transl.:
Theor.Math.Phys. 98 (1994) 147-161].
\item"[12]" Maltsev A.I., On a representation of nonassociative rings
[in Russian]. Uspekhi Matem. Nauk  7(1) (1952) 181-185 [reprinted in
``Selected Papers of Acad.A.I.Maltsev. I. Classical Algebra'', pp.328-331,
Moscow, Nauka, 1976].
\item"[13]" Drinfeld V.G., Quasi--Hopf algebras. Leningrad Math.J. 1(6)
(1990) 1419-1457; On quasitriangular quasi--Hopf algebras and a group closely
connected with $\Gal(\bar Q/Q)$. ibid. 2(4) (1991) 829-860;
Kosmann--Schwarzbach Y., Quasi--big\`ebres de Lie et groupes de Lie
quasi--Poisson. C.R.Acad.Sci.Paris, S\'erie I, 312 (1991) 391-394;
Jacobian quasi--bialgebras and quasi--Poisson Lie groups. Contemp.Math.132,
Amer.Math.Soc., Providence, 1992, pp.459-489; Bangoura M., Kosmann--Schwarzbach
Y., The double of a Jacobian quasi--bialgebra. Lett.Math.Phys. 28 (1993)
13-29; Bangoura M., Quasi--groupes de Lie--Poisson. C.R.Acad.Sci.Paris,
S\`erie I, 319 (1994) 975-978.
\item"[14]" Karasev M.V., Maslov V.P., Nonlinear Poisson brackets. Geometry
and quantization. Amer.Math.Soc., Providence, RI, 1993.
\item"[15]" Kobayashi Sh., Nomizu K., Foundations of differential geometry.
Interscience Publishers, 1963.
\item"[16]" Taylor M., Pseudo--differential operators. B., 1974; Duistermaat J.,
Fourier integral operators. N.Y., 1973; Treves F., Introduction to
pseudo--differential and Fourier integral operators. N.Y., 1980. 
\item"[17]" Zhelobenko D.P., Representations of reductive Lie algebras.
Moscow, Nauka, 1994.
\item"[18]" Ro\v cek M., Representation theory of the nonlinear $\SU(2)$
algebra. Phys. Lett.B 255 (1991) 554-557.
\item"[19]" Juriev D., Infinite dimensional geometry and quantum field theory
of strings. II. Infinite--dimensional noncommutative geometry of a 
self--interacting string field. Russian J.Math. Phys. 4(3) (1996) 287-314.
\item"[20]" Dixmier J., Alg\`ebres enveloppantes. Paris/Bruxelles/Montr\'eal,
Gauthier-Villars, 1974.
\item"[21]" Zamolodchikov A.B., Infinite additional symmetries in
two--dimensional quantum field theory [in Russian]. Teor.Matem.Fiz. 65(3)
(1986) 347-359.
\item"[22]" Juriev D., Algebraic structures of quantum projective field theory
related to fusion and braiding. Hidden additive weight. J.Math.Phys. 35 (1994)
3368-3379.
\item"[23]" Juriev D., Infinite dimensional geometry and quantum field theory 
of strings. III. Infinite--dimensional W--differential geometry of second 
quantized free string. J.Geom.Phys. 16 (1995) 275-300.
\item"[24]" Kac V.G., Infinite dimensional Lie algebras. Cambridge, Cambridge
Univ. Press, 1990; Feigin B.L., Fuchs D.B., Representations of the Virasoro 
algebra. In ``Representations of infinite dimensional Lie algebras''. Gordon
and Breach, 1991.
\endroster
P.S. The reader interested in the further elaboration of the ideas of the
articles [4,5] (and also of [11]) should be addressed to the report:
\roster
\item"[25]" Juriev D., On the description of a class of physical interactive
information systems [in Russian]. Report RCMPI/96-05 [e-version: 
mp\_arc/96-459].
\endRefs
\enddocument